\documentclass[12pt]{article}
\usepackage{amssymb,amsmath}

\newcommand{\be}{\begin{equation}} 
\newcommand{\ee}{\end{equation}}
\newcommand{\beq}{\begin{eqnarray}}
\newcommand{\eeq}{\end{eqnarray}}
\newcommand{\bl}{\begin{lemma}}
\newcommand{\el}{\end{lemma}}
\newcommand{\bc}{\begin{corollary}}
\newcommand{\ec}{\end{corollary}}
\newcommand{\bp}{\begin{prop}}
\newcommand{\ep}{\end{prop}}
\newcommand{\ba}{\begin{array}}
\newcommand{\ea}{\end{array}}

\newcommand{\wt}{\widetilde}
\newcommand{\wh}{\widehat}
\renewcommand{\mod}{\mathrm{mod}\;}
\newcommand{\tr}{\mathrm{tr}\;}
\newcommand{\la}{\label}
\newcommand{\ci}{\cite}

\newtheorem{theorem}{THEOREM}
\newtheorem{lemma}[theorem]{LEMMA}
\newtheorem{corollary}[theorem]{COROLLARY}
\newtheorem{prop}[theorem]{PROPOSITION}

\newcommand{\al}{\alpha}
\newcommand{\de}{\delta}
\newcommand{\De}{\Delta}

\newcommand{\ga}{\gamma}
\newcommand{\Ga}{\Gamma}

\newcommand{\si}{\sigma}

\newcommand{\lb}{\lambda}

\renewcommand{\th}{\theta}

\newcommand{\ve}{\varepsilon}

\newcommand{\bi}{\bibitem}

\newfont{\msbm}{msbm10 scaled\magstep1}
\newfont{\msbms}{msbm7 scaled\magstep1} 

\newcommand{\bbr}{\mbox{$\mbox{\msbm R}$}}

\newcommand{\bbz}{\mbox{$\mbox{\msbm Z}$}}

\renewcommand{\ep}{\varepsilon}

\begin{document}
\bigskip\bigskip\bigskip
\begin{center}
{\Large\bf Central spectral gaps of the almost Mathieu operator}\\
\bigskip\bigskip\bigskip
I. Krasovsky\\
\bigskip 
Department of Mathematics, Imperial College,
London, SW7 2AZ, UK

\end{center}
\bigskip\bigskip\bigskip

\noindent
{\bf Abstract.}
We consider the spectrum of the almost Mathieu operator $H_\al$ with frequency $\al$ and
in the case of the critical coupling. Let an irrational $\al$ be such that $|\al-p_n/q_n|<c q_n^{-\varkappa}$,
where $p_n/q_n$, $n=1,2,\dots$ are the convergents to $\al$, and
$c$, $\varkappa$ are positive absolute constants, $\varkappa<56$.
Assuming certain conditions on the parity of the coefficients of the continued fraction of $\al$,
we show that the central gaps of $H_{p_n/q_n}$, $n=1,2,\dots$, are inherited as spectral gaps 
of $H_\al$ of length at least $c'q_n^{-\varkappa/2}$, $c'>0$. 

\section{Introduction}
Let $H_{\al,\th}$ with $\al,\th\in(0,1]$ be the self-adjoint operator 
acting on $l^2(\bbz)$ as follows:
\be
(H_{\al,\th}\phi)(n)=\phi(n-1)+\phi(n+1)+2\cos 2\pi(\al n+\th)\phi(n),
\qquad n=\dots,-1,0,1,\dots
\ee
This operator is known as the almost Mathieu, Harper, or Azbel-Hofstadter operator.
It is a one-dimensional discrete periodic (for $\al$ rational) or 
quasiperiodic (for $\al$ irrational) Schr\"odinger operator which models 
an electron on the 2-dimensional square lattice in a perpendicular magnetic field. 
Analysis of the spectrum of $H_{\al,\th}$ (and its natural generalization
when the prefactor 2 of cosine, the coupling,  is replaced by an arbitrary real number $\lb$)
has been a subject of many investigations. In the present paper, we are concerned with
the structure of the spectrum of $H_{\al,\th}$ as a set. Denote by $a_j\in\bbz_+$, $j=1,2,\dots$,
the coefficients of the continued fraction of $\al$:
\[
\al=[a_1,a_2,\dots]=\frac{1}{a_1+\frac{1}{a_2+\cdots}}. 
\]
If $\al=p/q$ is rational, where $p$, $q$ are coprime, i.e. $(p,q)=1$, positive integers, 
there exists $n$ such that
\[
p/q=[a_1,a_2,\dots,a_n]=\frac{1}{a_1+\frac{1}{a_2+\cdots+\frac{1}{a_n}}}. 
\]

We denote by $S(\al)$ the union
of the spectra of  $H_{\al,\th}$ over all $\th\in (0,1]$ (Note, however, that if 
$\al$ is irrational, the spectrum of $H_{\al,\th}$ does not depend on $\th$).
If $\al=p/q$, $S(p/q)$ consists of $q$ bands separated by gaps.
As shown by van Mouche \ci{vM} and by Choi, Elliott, and Yui \ci{CEY},
all the gaps (with the exception of the centermost gap
when $q$ is even) are open. Much effort was expended to prove the conjectures
of \ci{Az,AA} that if $\al$ is irrational, the spectrum is a Cantor set.  
B\'ellissard and Simon proved in \ci{BS} that the spectrum of the generalized operator
mentioned above is a Cantor set
for an (unspecified) dense set of pairs $(\al,\lb)$ in $\bbr^2$. Helffer and Sj\"ostrand
\ci{HS}
proved the Cantor structure and provided an analysis of gaps
in the case when all the coefficients $a_j$'s of $\al$ are sufficiently large. 
Choi, Elliott, and Yui \ci{CEY} showed that in the case of $\al=p/q$, each open gap
is at least of width $8^{-q}$ (this bound was improved in \ci{AJ} 
to $e^{-\ep q}$ with any $\ep>0$ for $q$ sufficiently large) which, 
together with a continuity result implies that {\it all} admissible
gaps are open (in particular, the spectrum is a Cantor set)
if $\al$ is a Liouvillian number whose convergents $p/q$ satisfy $|\al-p/q|<e^{-C q}$.
Last \ci{L} showed that $S(\al)$ has Lebesgue measure zero (and hence,
since $S(\al)$ is closed and known not to contain isolated points, a Cantor set) for all
$\al=[a_1,a_2,\dots]$ such that the sequence $\{a_j\}_{j=1}^\infty$ is unbounded.
The set of such $\al$'s has full measure $1$. On the other hand, it was shown by Puig
\ci{P} that in the generalized case $\lb\neq \pm 2,0$, the spectrum is a Cantor set for
$\al$ satisfying a Diophantine condition.  
Finally, Avila and Krikorian \ci{AK} completed the proof that
the spectrum for $\lb=2$ has zero measure, and hence a Cantor set, 
for all irrational $\al$'s; moreover, the proof of 
the fact that the spectrum is a Cantor set for all real $\lb\neq 0$ and irrational $\al$
was completed by Avila and Jitomirskaya in \ci{AJ}. The measure of the spectrum for any irrational 
$\al$ and real $\lb$ is $|4-2|\lb||$: in the case $\lb\neq\pm2$,
proved for a.e. $\al$ also in \ci{L} and for all irrationals in \ci{JK}. Also available are bounds 
on the measure of the union of all gaps, see \ci{DS,L2,KK}. Furthermore, see \ci{L15}
for a recent work on the Hausdorff dimension of the spectrum, and \ci{LY}, on the question
of whether all admissible gaps are open.  


In order to have a quantitative description of the spectrum, one would like to know 
if the exponential $e^{-\ep q}$ estimates for the sizes of the individual gaps can be improved 
at least for some of the gaps.

In this paper we provide a power-law estimate $Cq^{-\kappa}$, $\kappa<28$, for the widths of
{\it central} gaps of $S(p/q)$, i.e. the gaps around the centermost band (Theorem 3 below),
on a parity condition for the coefficients $a_k$ in $p/q=[a_2,a_2,\dots,a_n]$.

From this result we deduce that $S(\al)$ has an infinite number of power-law bounded gaps
for any irrational $\al=[a_1,a_2,...]$ admitting a power-law approximation by its convergents 
$p_n/q_n=[a_1,a_2,\dots,a_n]$ and with a parity condition on $a_j$'s (Theorem 4 below).
These gaps are inherited from the central ones of $S(p_n/q_n)$, $n=1,2,\dots$.

First, let $\al=p/q$, $(p,q)=1$. A standard object used for the analysis of $H_{\al,\th}$
is the discriminant
\begin{align}
\si(E)=&-\tr\left\{
\begin{pmatrix} E-2\cos(2\pi p/q+\pi/2q)& -1\cr 1 & 0 \end{pmatrix}
 \begin{pmatrix} E-2\cos(2\pi 2p/q+\pi/2q)& -1\cr 1 & 0 \end{pmatrix}\cdots\right.\nonumber\\
&\left.
 \begin{pmatrix} E-2\cos(2\pi qp/q+\pi/2q)& -1\cr 1 & 0 \end{pmatrix}
\right\},\la{discrim}
\end{align}
a polynomial of degree $q$ in $E$ with the property that $S(p/q)$ is the
image of $[-4,4]$ under the inverse of the mapping $\si(E)$.
The fact that $S(p/q)$ consists of $q$ bands separated by $q-1$ open gaps (except for the centermost 
empty gap for $q$ even)
means that all the zeros of $\si(E)$ are simple, in all the maxima 
the value of $\si(E)$ is strictly larger than $4$, while in all the minima, strictly less 
than $-4$ (except for $E=0$ for $q$ even, where $|\si(0)|=4$ and the derivative $\si'(0)=0$).
Note an important fact that $\si(E)=(-1)^q\si(-E)$,
and hence $S(p/q)$ is symmetric w.r.t. $E=0$.

In what follows, we assume that $q$ is odd. The case
of even $q$ can be considered similarly.
Let us number the bands from left to right, from $j=-(q-1)/2$ to $j=(q-1)/2$.
Let $\lb_j$ denote the centers of the bands, i.e. $\si(\lb_j)=0$.
Note that, by the symmetry of $\si(E)$, $\lb_0=0$.
Let $\mu_j$ and $\eta_j$ denote the edges of the bands, i.e. $|\si(\mu_j)|=|\si(\eta_j)|=4$, 
assigned as follows. If $q=4k+3$, $k=0,1,\dots$, we set $\si(\mu_j)=4$, $\si(\eta_j)=-4$ for all $j$.
(In this case the derivative $\si'(0)>0$, as follows from the fact that $\si(E)<0$ for all $E$ sufficiently large.)
If $q=4k+1$, $k=0,1,\dots$, we set $\si(\mu_j)=-4$, $\si(\eta_j)=4$ for all $j$. 
(In this case the derivative $\si'(0)<0$.) Thus, in both cases, the bands are
$B_j=[\eta_j,\mu_j]$ for $|j|$ even, and $B_j=[\mu_j,\eta_j]$ for $|j|$ odd.

Let $w_j=\mu_j-\lb_j$, $w'_j=\lb_j-\eta_j$ for $|j|$ even, 
and $w_j=\eta_j-\lb_j$, $w'_j=\lb_j-\mu_j$ for $|j|$ odd.
Thus, the width of the $j$'s band is always $w_j+w'_j$. By the symmetry, for
the centermost band $B_0=[\eta_0,\mu_0]$, $w_0=w'_0$, and in general $w_j=w'_{-j}$.

For any real $\al$, denote the {\it gaps} of $S(\al)$ by $G_j(\al)$ and their length by $\De_j(\al)$.
For $\al=p/q$, we order them in the natural way, namely,
\begin{align}
G_j&=(\mu_j,\mu_{j+1}),\quad \De_j=\mu_{j+1}-\mu_j,\quad\mbox{for $|j|$ even},\\
G_j&=(\eta_j,\eta_{j+1}),\quad \De_j=\eta_{j+1}-\eta_j,\quad\mbox{for $|j|$ odd}.
\end{align}

By the symmetry, $\De_j=\De_{-j-1}$ for $ 0\le j<(q-1)/2$.

In Section 2, we prove

\noindent {\bf Lemma 1} (Comparison of widths for gaps and bands)
{\it Let $q\ge 3$ be odd.
There hold the inequalities
\begin{align}\la{L11}
\De_0>\left(\frac{w_0}{4}\right)^2,\qquad
\De_j>&\frac{w_j^2}{4C_0^{2(j+1)}},\qquad 1\le j<\frac{q-1}{2},\\
\De_j>&\left(\frac{w_0}{8}\right)^{2j},\qquad 1\le j<\frac{q-1}{2},\la{L12}
\end{align}
where $C_0=1+2e/(\sqrt{5}-1)=5.398\dots$.
}

\noindent{\bf Remark} The inequalities of Lemma 1 are better for small $j$, i.e.,
for central gaps and bands, which is the case we need below. 
For large $j$, note the following estimate which one can deduce using the technique
of Last \ci{L}: $\De_j>\min\{{w_j}^2,{w'_{j+1}}^2\}/(4q)$, $0\le j<(q-1)/2$.

The inequality (\ref{L12}) gives us a lower bound for the width of the $j$'s gap
provided an estimate for the width of the $0$'s band can be established.
Such an estimate is given by  

\noindent{\bf Lemma 2} (Bound for the width of the centermost band)
{\it Let $q\ge 1$, $p/q=p_n/q_n=[a_1,a_2,\dots,a_n]$, where $a_1$ is odd and $a_k$, 
$2\le k\le n$ are even.
Then there exist absolute constants $1<C_1<14$ and $1<C_2<e^{10}$ such that for the 
derivative of $\si(E)$ at zero
\be\label{si0}
|\si'(0)|< C_2 q^{C_1},
\ee
and half the width of the centermost band of $S(p/q)$
\be\label{w0}
w_0 \ge {4\over |\si'(0)|} > 4C_2^{-1} q^{-C_1}.
\ee

If, in addition, $q_{k+1}\ge q_k^{\nu}$, for some $\nu>1$ and all $1\le k\le n-1$, then for any $\ve>0$
there exists $Q=Q(\ep,\nu)$ such that if $q>Q$, 
\be\la{siw-new}
|\si'(0)|< q^{5+\ga_0+\ve},\qquad w_0 > 4q^{-(5+\ga_0+\ve)},
\ee
where $\ga_0$ is Euler's constant.
}

\medskip
\noindent{\bf Remark} 
The bounds on $C_1$, $C_2$ can be somewhat improved.

This lemma is proved in Section 3.
The inequalities (\ref{L11}), (\ref{L12}), and especially (\ref{si0}) are the main technical results of this paper.

Combination of Lemmata 1 and 2  immediately yields

\noindent{\bf Theorem 3} (Bound for the widths of the gaps) {\it
Let $q\ge 3$, $p/q=[a_1,a_2,\dots,a_n]$, where $a_1$ is odd and $a_k$, $2\le k\le n$, are even.
Then, with $C_k$, $k=1,2$, from Lemma 2,
the width of the $j$'s gap of $S(p/q)$ is
\be
\De_0>\left(\frac{1}{C_2 q^{C_1}}\right)^2,\qquad
 \De_j >\left(\frac{1}{2C_2 q^{C_1}}\right)^{2j},\qquad 1\le j < \frac{q-1}{2}.
\ee
}

\noindent{\bf Remark} The improvements for large $q$ on the additional condition $q_{k+1}>Cq_k^{\nu}$
are obvious from (\ref{siw-new}).

A consequence of this is the following theorem proved in Section 4.

\noindent{\bf Theorem 4} {\it
There exists an absolute $C_3>0$ such that the following holds.
Let $\al=[a_1,a_2,\dots]\in (0,1)$ be an irrational such that $a_1$ is odd, $a_k$, $k\ge 2$, are even,
and such that 
\be\la{alphacond}
\left|\al-\frac{p_n}{q_n}\right|<\frac{1}{C_3 q_n^\varkappa}, \qquad \varkappa=4C_1, 
\ee
for all $p_n/q_n=[a_1,a_2,\dots,a_n]$, $n=1,2,\dots$, where $C_1$ is the constant from Lemma 2. 

Then

\noindent
(a) The interior of the centermost band $B_0$ of $S(p_n/q_n)$ contains the centermost band
and the closures of the gaps $G_0$, $G_{-1}$ of $S(p_{n+1}/q_{n+1})$, $n=1,2,\dots$

\noindent
(b) There exist distinct gaps $G_{n,j}(\al)$, $n=1,2,\dots$, $j=1,2$, of $S(\al)$,
such that the intersections $G_{n,1}(\al)\cap G_{-1}(p_n/q_n)$, $G_{n,2}(\al)\cap G_{0}(p_n/q_n)$,
$n=1,2,\dots$ are non-empty
and the length of the gap $G_{n,j}(\al)$ 
\be\la{Dewidth}
\De_{n,j}(\al)=|G_{n,j}(\al)|
\ge |G_{n,j}(\al)\cap G_{j-2}(p_n/q_n)|
>\frac{1}{C_4 q_n^{\varkappa/2}},\qquad n=1,2,\dots,\quad j=1,2,
\ee
for some absolute $C_4 >0$, where $|A|$ denotes the Lebesgue measure of $A$.

\noindent
(c) Let $\ve>0$,  replace $C_3$ by $2$, and set $\varkappa=4(5+\ga_0+\ve)$
in (\ref{alphacond}). Then there exists $n_0=n_0(\ve)$ such that (a) and (b) hold for all $n=n_0,n_0+1,\dots$
(instead of $n=1,2,\dots$) with $C_4$ replaced by $2$, and with $\varkappa/2$ in 
(\ref{Dewidth}) replaced by $2(5+\ga_0)+\ve$.
}

\medskip
\noindent{\bf Remarks} 

\noindent
1) The statements (a), (b) of the theorem hold a fortiori for $\varkappa=4\cdot 14=56$ and for any larger
$\varkappa$.
It is easy to provide explicit examples of irrationals satisfying the conditions of Theorem 4:
take $\varkappa=56$, any odd $a_1$, and even $a_{n+1}$ such that 
$a_{n+1}>C_3 q_n^{\varkappa-2}$, $n\ge 1$. Indeed, in this case,
\[
\left|\al-\frac{p_n}{q_n}\right|<\frac{1}{q_{n}q_{n+1}}<\frac{1}{a_{n+1}q_n^2}<
\frac{1}{C_3 q_n^\varkappa}.
\]

\noindent
2)  Note that the parity condition on $a_j$'s implies, in particular, that all $q_n$'s are odd.
This condition can be relaxed in all our statements. For example, we can allow
a finite number of $a_j$'s to be odd at the expense of excluding some $G(p_n/q_n)$'s
from the statement of Theorem 4 and worsening the bound
on $C_1$. Note that in Lemma 2 we need $q$ to be odd in order to use
the estimate (\ref{si0}) on $\si'(0)$ to obtain (\ref{w0}).
One could obtain a bound on $w_0$ for even $q$ by providing an estimate
on the second derivative $\si''(0)$ in this case: for $q$ even $\si'(0)=0$. 
The parity condition we assume in this paper allows the best estimates and simplest proofs.

\noindent
3) In Theorem 4, we only use Theorem 3 for $j=0$, i.e., for the 2 centermost gaps.
One can extend the result of Theorem 4, with appropriate changes, to more than 2 (at least a finite number) of central gaps of $S(p_n/q_n)$.

\noindent
4) We can take $C_3=4^2 60^2 C_2^4$, $C_4=2C_2^2$, in terms of the constant $C_2$
from Lemma 2. 

\noindent
5) The statement (a) of the theorem holds already for $\varkappa=2C_1$.

\section{Proof of Lemma 1}
Assume that $q=4k+3$, $k=0,1,\dots$. (A proof in the case $q=4k+1$ is almost identical.)
Let
\[
s={q-1\over 2}.
\]
In our notation, we can write
\be
\si(E)=\prod_{k=-s}^s (E-\lb_k),\qquad
\si(E)-4=\prod_{k=-s}^s(E-\mu_k),\qquad
\si(E)+4=\prod_{k=-s}^s(E-\eta_k).
\ee
Setting in the last 2 equations $E=\lb_j$, we obtain the useful identities
\be\la{prodid}
4=\prod_{k=-s}^s|\lb_j-\mu_k|,\qquad 4=\prod_{k=-s}^s|\lb_j-\eta_k|,\qquad 
-s\le j\le s.
\ee

Fix $0\le j\le s$
(by the symmetry of the spectrum, it is sufficient to consider only nonnegative $j$).
It was shown by Choi, Elliott, and Yui \ci{CEY} that
\be\la{Ell}
\prod_{k\neq j}|\mu_j-\mu_k|\ge 1,\qquad  \prod_{k\neq j}|\eta_j-\eta_k|\ge 1.
\ee

For simplicity of notation, we assume from now on that $j<s-1$: the extension to $j=s-1$ is obvious. 
Let $j\ge 0$ be even. By the first inequality in (\ref{Ell}), we can write
\begin{align}\la{prodint}
1\le|\si'(\mu_j)|&=\prod_{k\neq j}|\mu_j-\mu_k|=\nonumber\\
&|\mu_j-\mu_{j+1}|\frac{\prod_{k=-s}^s|\lb_j-\mu_k|}{|\lb_j-\mu_j||\lb_j-\mu_{j+1}|}
\prod_{k=-s}^{j-1}\left|1+\frac{\mu_j-\lb_j}{\lb_j-\mu_k}\right| 
\prod_{k=j+2}^{s}\left|1-\frac{\mu_j-\lb_j}{\mu_k-\lb_j}\right|.
\end{align}
According to our notation, $\mu_{j+1}-\mu_j=\De_j$, $\mu_j-\lb_j=w_j$,
$\mu_{j+1}-\lb_j=w_j+\De_j$. Recalling the first identity in (\ref{prodid}) and 
rearranging the last product in (\ref{prodint}), we continue (\ref{prodint})
as follows
\be\la{prodint2}
=\frac{4\De_j}{w_j(w_j+\De_j)}
\frac{\prod_{k=-s}^{j-1}\left|1+\frac{w_j}{\lb_j-\mu_k}\right|}
{\prod_{k=j+2}^{s}\left|1+\frac{w_j}{\mu_k-\mu_j}\right|}<
\frac{4\De_j}{w_j(w_j+\De_j)}
\frac{\prod_{k=-s}^{j-1}\left|1+\frac{w_j}{\lb_j-\mu_k}\right|}
{\prod_{k=j+2}^{s}\left|1+\frac{w_j}{\mu_k-\lb_j}\right|},
\ee
because $\mu_j>\lb_j$. 

Now note that, by the symmetry of the spectrum,
\be
|\mu_{j+\ell}-\lb_j|<|\mu_{-j-\ell-1}-\lb_j|,\qquad \ell=2,3,\dots
\ee
Therefore, the r.h.s. of (\ref{prodint2}) is
\be\la{p3}
< \frac{4\De_j}{w_j(w_j+\De_j)}\prod_{k=-j-2}^{j-1}\left|1+\frac{w_j}{\lb_j-\mu_k}
\right|\frac{1}{1+{w_j\over \mu_s-\lb_j}}.
\ee

In the case $j=0$, we now use the symmetry
\be
w_0=w_0'<\lb_0-\mu_{-k},\qquad k\ge 1,
\ee
to obtain from (\ref{p3})
\[
1<\frac{16\De_0}{w_0(w_0+\De_0)}<\frac{16\De_0}{w_0^2},
\]
which gives the first inequality in (\ref{L11}).
 
 In general, however, we need to compare $w_j$ and $w'_j$ to estimate (\ref{p3}). 
According to equations (3.11), (3.12) of Last \ci{L},
\be\la{Last-eq1}
w_j, w'_j < e\ell_j,\qquad \ell_j=\frac{4}{|\si'(\lb_j)|},
\ee
and further, by equations (3.27), (3.28) of \ci{L},
\be\la{Last-eq2}
\frac{\sqrt{5}-1}{2}\ell_j<w_j, w'_j. 
\ee
(In fact, more is shown in \ci{L}: for each pair of widths $w_j$, $w'_j$, at least one of them is larger
than $\ell_j$.)  

Therefore,
\be
w_j<c_1 w_j',\qquad c_1=\frac{2e}{\sqrt{5}-1},\qquad 0\le j\le s.
\ee
Furthermore, it is obvious that 
\be
\frac{w_j'}{\lb_j-\mu_k}<1,\qquad k=-j-2,\dots,j-1.
\ee
Therefore, we have for the product in (\ref{p3}):
\be
\prod_{k=-j-2}^{j-1}\left|1+\frac{w_j}{\lb_j-\mu_k}\right|<(1+c_1)^{2(j+1)},
\ee
and since $\mu_s-\lb_j>0$, (\ref{p3}) finally gives
\be
1<\frac{4\De_j}{w_j(w_j+\De_j)}(1+c_1)^{2(j+1)},
\ee
from which the inequality (\ref{L11}) with even $j$
easily follows. 

\noindent
{\bf Remark} Last's equation (\ref{Last-eq1}) together with the Last-Wilkinson formula \ci{L,LW}
\be\la{LW}
\sum_{j=-s}^s |\si'(\lb_j)|^{-1}=1/q
\ee
implies \ci{L}  that the measure of the spectrum $S(p/q)$ is at most $8e/q$ and that for any $j$, 
\be\la{w-est}
w_j<4e/q.
\ee


Now consider $j$ odd, $0<j<s$. Using the second inequalities
in (\ref{Ell}) and (\ref{prodid}), we obtain similarly to (\ref{prodint2}),
\be
1< \frac{4\De_j}{w_j(w_j+\De_j)}
\frac{\prod_{k=-s}^{j-1}\left|1+\frac{w_j}{\lb_j-\eta_k}\right|}
{\prod_{k=j+2}^{s}\left|1+\frac{w_j}{\eta_k-\lb_j}\right|},
\ee
and since
\be
|\eta_{j+\ell}-\lb_j|<|\eta_{-j-\ell-1}-\lb_j|,\qquad \ell=2,3,\dots,
\ee
we obtain the inequality (\ref{L11}) for $j$ odd in a similar way.

Let again $j$ be even, $0< j<s-1$.
In order to compare $\De_j$ with the width of the centermost band and, thus, obtain 
(\ref{L12}), we write instead of (\ref{prodint}) the following:
\begin{align}\la{bprodint}
1\le|\si'(\mu_j)|&=\prod_{k\neq j}|\mu_j-\mu_k|=\\
&|\mu_j-\mu_{j+1}|\frac{\prod_{k=-s}^s|\lb_0-\mu_k|}{|\lb_0-\mu_j||\lb_0-\mu_{j+1}|}
\prod_{k=-s}^{j-1}\left|1+\frac{\mu_j-\lb_0}{\lb_0-\mu_k}\right| 
\prod_{k=j+2}^{s}\left|1-\frac{\mu_j-\lb_0}{\mu_k-\lb_0}\right|.
\end{align}
Proceeding in a similar way as before, and using the inequalities
\be
|\mu_{j+\ell}-\lb_0|<|\mu_{-j-\ell-1}-\lb_0|,\qquad \ell=2,3,\dots,
\ee
we obtain
\begin{align}
1&< \frac{4\De_j}{|\lb_0-\mu_j||\lb_0-\mu_{j+1}|}
\prod_{k=-j-2}^{j-1}\left|1+\frac{\mu_j-\lb_0}{\lb_0-\mu_k}\right|
<\frac{16\De_j}{|\lb_0-\mu_j||\lb_0-\mu_{j+1}|}
\prod_{k=-j}^{j-1}\left|1+\frac{\mu_j-\lb_0}{\lb_0-\mu_k}\right|\nonumber\\
&=\frac{16\De_j}{|\lb_0-\mu_j||\lb_0-\mu_{j+1}|}
\left|\frac{\mu_j-\mu_0}{\lb_0-\mu_0}\right|\left|\frac{\mu_j-\mu_1}{\lb_0-\mu_1}\right|
\prod_{{k=-j\atop k\neq 0,1}}^{j-1}\left|\frac{\mu_j-\mu_k}{\lb_0-\mu_k}\right|\nonumber\\
&<\frac{16\De_j}{|\lb_0-\mu_0||\lb_0-\mu_{1}|}
\prod_{{k=-j\atop k\neq 0,1}}^{j-1}\left|\frac{\mu_j-\mu_k}{\lb_0-\mu_k}\right|,
\end{align}
and since (note that $S(\al)\in[-4,4]$)
\[
\left|\frac{\mu_j-\mu_k}{\lb_0-\mu_k}\right|<{8\over w_0},
\]
we obtain
\be
1< \frac{\De_j}{4}\left({8\over w_0}\right)^{2j},
\ee
which gives an inequality slightly better than (\ref{L12}) for $j$ even.
Finally, we establish (\ref{L12}) for $j$ odd by starting (instead of (\ref{bprodint}))
with the inequality
$1\le|\si'(\eta_{j+1})|=\prod_{k\neq j+1}|\eta_{j+1}-\eta_k|$ and arguing similarly. 

\noindent{\bf Remark} Using the Last estimate (\ref{Last-eq1})
\be
{4e\over w_j}> |\si'(\lb_j)|=\prod_{k\neq j}|\lb_j-\lb_k|,
\ee
one can establish, in a way similar to the argument above, inequalities of the type
\be
\De_j+w'_{j+1}>\frac{w_j}{C^j},
\ee
with some absolute constant $C>0$. $\Box$

\section{Proof of Lemma 2}
As noted in a remark following Theorem 4, the parity conditions imposed on $p/q$ in Lemma 2
imply, in particular, that $q$ is odd. 
It follows from the symmetry of the discriminant $\si(E)=-\si(-E)$ in this case 
that the maximum of the absolute value
of the derivative $\si'(E)$ in the $j=0$ band is at $E=0$. Therefore,
\be\la{w-discr}
w_0 \ge {4\over |\si'(0)|}
\ee
(with the equality only for $q=1$),
and hence, in order to prove Lemma 2, it remains to obtain the inequality (\ref{si0}).

If $q=1$, we have $\si(E)=-E$, and the result is trivial. Assume now that $q$ is any 
(even or odd) integer larger than $1$.
We start with the following representation of $\si(E)$ in terms of a $q\times q$ Jacobi matrix
with the zero main diagonal:
\be\la{whH}
\si(E)=\det(\wh H-EI),
\ee
where $I$ is the identity matrix, and $\wh H$ is a $q\times q$ matrix $\wh H_{j\,k}$, 
$j,k=1,\dots q$, where
\be
\wh H_{j\, j+1}=\wh H_{j+1\, j}=2\sin\left(\pi\frac{p}{q}j\right),\qquad j=1,\dots,q-1,
\ee
and the rest of the matrix elements are zero. For a proof, see e.g. the appendix of \ci{Kprb}.
(This is related to a matrix representation for the almost Mathieu operator
corresponding to the chiral gauge of the magnetic field potential, noticed 
by several authors \ci{MZ,KH0d,WZmpl}.) 
The absence of the main diagonal in $\wh H$ allows us to obtain a simple expression
for the derivative $\si'(E)$ at $E=0$. If $q$ is even, it is easily seen that $\si'(E)=0$.
If $q$ is odd, we denote $s=(q-1)/2$ and 
immediately obtain from (\ref{whH}) 
(henceforth we set $\prod_{j=a}^b\equiv 1$ and $\sum_{j=a}^b\equiv 0$ if $a>b$):
\be\label{si0exact}
\si'(0)=(-1)^{s}\sum_{k=0}^{s}
\left[\prod_{j=1}^k 2\sin{\pi p\over q}(2j-1)\prod_{j=k+1}^{s}2\sin{\pi p\over q}2j
\right]^2.
\ee

From now on, we assume that $q\ge 3$ is odd unless stated otherwise.

\noindent
{\bf Remark} Using the identity $\prod_{j=1}^{(q-1)/2}2\sin{\pi p\over q}2j=q$, we can represent 
(\ref{si0exact}) in the form 
\be
\si'(0)=(-1)^{s}q\left(1+\sum_{k=1}^{s}\prod_{j=1}^k
\frac{\sin^2{\pi p\over q}(2j-1)}{\sin^2{\pi p\over q}2j}\right),
\ee
which exhibits the fact that $|\si'(0)|>q$. This is in accordance with the Last-Wilkinson formula
(\ref{LW}).

\medskip

Thus we have
\be\la{abs-si}
|\si'(0)|=\sum_{k=0}^{s}\exp\{L_k\},
\ee
where
\begin{align}
\frac{1}{2}L_k=&\sum_{j=1}^k\ln\left|2\sin{\pi p\over q}(2j-1)\right|+
\sum_{j=k+1}^{s}\ln\left|2\sin{\pi p\over q}2j\right|\la{Lm1}\\
=&\sum_{j=-s+k+1}^k\ln\left|2\sin{\pi p\over q}(2j-1)\right|=
\sum_{j=1}^s\ln\left|2\sin{\pi p\over q}2(j+k)\right|.\la{L}
\end{align}
Here we changed the summation variable $j=s+j'$ in the second sum in (\ref{Lm1})
to obtain the first equation in (\ref{L}), and then changed the variable $j=k-s+j'$ to obtain
the final equation in (\ref{L}).

We will now analyze $L_k$. Using the Fourier expansion, we can write
\be
L_k=
- 2\sum_{j=1}^s\sum_{n=1}^\infty {1\over n}\cos 4n {\pi p\over q}(j+k).
\ee
Representing $n$ in the form $n=qm+\ell$, where $\ell=1,2,\dots,q-1$ for $m=0$, and
$\ell=0,1,\dots,q-1$ for $m=1,2,\dots$, we have
\be\la{Lmed}
L_k=
-\sum_{m=1}^\infty\left({q-1\over qm}-{1\over q}\sum_{\ell=1}^{q-1}{1\over m+\ell/q}
F(\ell,k)\right)+S_k,
\ee
where
\be\la{Sk}
S_k=\sum_{\ell=1}^{q-1}{1\over\ell}F(\ell,k),
\ee
and
\be\la{F}
F(\ell,k)=-2\sum_{j=1}^s \cos  4\ell {\pi p\over q}(j+k)=
-2\Re\sum_{j=1}^s \exp\{4\pi i {p\over q}(j+k)\ell\}=
\frac{\cos \pi{p\over q}\ell(4k+1)}{\cos \pi{p\over q}\ell}.
\ee
Note that since $1<j+k\le q-1$, $(p,q)=1$, and $q$ is odd, we have that
$(2j+2k)p\not\equiv 0(\mod q)$.

For $0\le k\le s$,
\be\la{Fsum}
\sum_{\ell=1}^{q-1}F(\ell,k)=
-2\Re\left(\sum_{j=1}^s \sum_{\ell=1}^{q-1} \exp\{4\pi i {p\over q}(j+k)\ell\}
+1-1\right)=q-1.
\ee
We will use this fact later on.

Recall that 
\be\la{e1}
\sum_{m=1}^M{1\over m}=\ln M+\ga_0+o(1),\qquad M\to\infty,
\ee
where $\ga_0=0.5772\dots$ is Euler's constant. Recall also Euler's $\psi$-function
\be\la{e2}
\psi(x+1)=\frac{\Gamma'(x+1)}{\Gamma(x+1)}=
-\ga_0-\sum_{m=1}^\infty\left({1\over m+x}-{1\over m}\right),\qquad x\ge 0.
\ee
The function $\psi(x)$ continues to a meromorphic function in the complex plane with
first-order poles at nonpositive integers $x=0,-1,-2,\dots$. The function $\psi(x)$
satisfies the equation
\[
\psi(x+1)=\psi(x)+{1\over x}.
\]

Expressions (\ref{e1}), (\ref{e2}) imply
\be\la{e3}
 \sum_{m=1}^M{1\over m+x}=\sum_{m=1}^M\left({1\over m+x}-{1\over m}\right)+
\sum_{m=1}^M{1\over m}=\ln M-\psi(x+1)+o(1),\qquad  M\to\infty,
\ee
uniformly for $x\in[0,1]$, in particular.
We now rewrite (\ref{Lmed}) in the form
\be
L_k=-\lim_{M\to\infty}
\left[{q-1\over q}\sum_{m=1}^M {1\over m} -{1\over q}\sum_{\ell=1}^{q-1}F(\ell,k)
\sum_{m=1}^M{1\over m+\ell/q}\right]+S_k.
\ee
Substituting here (\ref{e1}), (\ref{e3}), and using (\ref{Fsum}), we finally obtain
\be\la{Lend}
L_k=- {q-1\over q}\ga_0-
{1\over q}\sum_{\ell=1}^{q-1}F(\ell,k)\psi(1+\ell/q)+S_k.
\ee

We will now provide an upper bound for the absolute values of sums in this expression. 
First, note that the derivative $\psi'(x)\ge 0$, $x\in[1,2]$, and
$\psi(1)=-\ga_0$, $\psi(2)=\psi(1)+1=1-\ga_0$. Therefore,
\[
\max_{x\in[1,2]}|\psi(x)|=\max\{|\psi(1)|,|\psi(2)|\}=\ga_0.
\]
Thus, for the first sum in the r.h.s. of (\ref{Lend}) we have
\begin{align}
\left|{1\over q}\sum_{\ell=1}^{q-1}F(\ell,k)\psi(1+\ell/q)\right|&\le
{\ga_0\over q}\sum_{\ell=1}^{q-1}|F(\ell,k)|=
{\ga_0\over q}\sum_{\ell=1}^{q-1}{1\over |\cos(\pi p\ell/q)|}\nonumber\\
&=\frac{4\ga_0}{q}\sum_{m=0}^{s-1}\frac{1}{|1-e^{i(2m+1)\pi/q}|}
< \ga_0\sum_{m=0}^{s-1}\frac{1}{m+1/2}\nonumber\\
&<\ga_0\left(2+\int_0^{s-1}\frac{dx}{x+1/2}\right)<
\ga_0\left(\ln q+2\right),\qquad q\ge3.\la{1sumest}
\end{align}

We will need a more subtle estimate for 
\be\la{Sk2}
S_k=
\sum_{m=1}^{q-1}\frac{1}{m}F(m,k)=\sum_{m=1}^{q-1}\frac{1}{m}
\frac{\cos \pi{p\over q}m(4k+1)}{\cos \pi{p\over q}m}
\ee
in the r.h.s. of (\ref{Lend}). We follow a method of Hardy and Littlewood \ci{HL1930} (see also 
\ci{Spencer1939}).
It relies on a recursive application of a suitably constructed contour integral.

For $q\ge 3$ odd, $(p,q)=1$, let
\be\la{1int}
I(p/q,\ga)= -2\int_{\Ga_{q}}\frac{e^{(1+p/q)z}}{(1+e^{z p/q})(1-e^z)}\frac{e^{-\ga z}}{z}dz,
\qquad \frac{1}{2}\frac{p}{q}\le\ga\le1+\frac{1}{2}\frac{p}{q},
\ee
where the contour $\Ga_{q}$ are the 2 direct lines parallel to the real axis given by:
(1) $\pi i/2+x$, $x\in\bbr$, oriented from $-\infty$ to $+\infty$;
(2) $2\pi i(q-1/4)+x$,  $x\in\bbr$, oriented from $+\infty$ to $-\infty$.
Note that the choice of $\ga$ in (\ref{1int}) ensures that the integral converges
both at $+\infty$ and $-\infty$.

Now again for $q\ge 3$ odd, , $(p,q)=1$, let
\be\la{SumS}
S(p/q,\ga)=\sum_{m=1}^{q-1}\frac{e^{\pi i\frac{p}{q}m-2\pi i\ga m}}{m\cos{(\pi{p\over q}m)}}.
\ee
Note that the sum (\ref{Sk2})
\be
S_k=\Re S(p/q,-2kp/q),
\ee
and that the denominators in (\ref{SumS}) are nonzero.

We will also need the following auxiliary sum, $(p,q)=1$,
\be\la{SumT}
T(p/q,\ga,\de)=2\sum_{n=1}^q\frac{(-1)^{\de}e^{2\pi i{p\over q}(n-{1\over 2})}}{1-(-1)^{\de}
e^{2\pi i{p\over q}(n-{1\over 2})}}\frac{e^{-2\pi i\ga(n-{1\over 2})}}{n-{1\over 2}},\qquad \de=0,1,
\ee
where we assume that $p$ is odd if $\de=0$ and that
$p$ and $q$ have opposite parities if $\de=1$.
These conditions imply that $p(2n-1)\neq q(2m-\de)$, $m,n\in\bbz$, and therefore the denominators in
(\ref{SumT}) are nonzero.

With this notation we have

\medskip
\noindent
{\bf Lemma 5}. {\it Let $q$ be odd, $(p,q)=1$, $p>0$, $q>1$. 
Let $p/q$ have the following continued fraction:
\be\la{cf}
\frac{p}{q}=\frac{1}{a+\frac{p'}{q'}},\qquad p'\ge 0, \quad q'>p'.
\ee
Then 
\be\la{1intrec}
I(p/q,\ga)=S(p/q,\ga)-(-1)^{\ep'}T(p'/q',\ga',a\,\mod2),\qquad \ga'=\frac{q}{p}\ga (\mod 1),
\ee
where $\ep'=0$ if $\ga'=\frac{q}{p}\ga (\mod 2)$, and $\ep'=1$ otherwise.

Moreover, there holds the bound
\be\la{1intest}
|I(p/q,\ga)|<4\ln\frac{q}{p}+\frac{5}{e\pi p}+
\beta,\qquad \beta=4(e^{-1}+\mathrm{arcsinh}(4/\pi))=5.719\dots
\ee
}

\noindent
{\bf Remarks} 

\noindent
1) One can take any $\ga'$ satisfying the congruence in (\ref{1intrec}).

\noindent
2) The bound in (\ref{1intest}) can be somewhat decreased by improving (\ref{esti1}), (\ref{1intest1}) below.
Similar can be achieved in (\ref{2intest}) below.

\noindent
{\bf Proof}
Consider the integral $I_{\Gamma}(p/q,\ga)$, which has the same 
integrand as in (\ref{1int}), but with integration over some contour $\Gamma$.
Let $\Gamma_{q,\xi}$ be the following quadrangle traversed in the positive direction:
$\Gamma_{q,\xi}=\cup_{j=1}^4 \Gamma^{(j)}_{q,\xi}$, $\xi>0$, 
where 
$\Gamma^{(1)}_{q,\xi}=[\pi i/2-\xi, \pi i/2+\xi]$,
$\Gamma^{(2)}_{q,\xi}=[\pi i/2+\xi,2\pi i(q-1/4)+\xi]$,
$\Gamma^{(3)}_{q,\xi}=[2\pi i(q-1/4)-\xi,2\pi i(q-1/4)+\xi]$,
$\Gamma^{(4)}_{q,\xi}=[2\pi i(q-1/4)-\xi,\pi i/2-\xi]$.

Recalling the conditions on $\ga$ in (\ref{1int}), we
first note that on the vertical segments, for some constants $C$ which may depend on $p$, $q$,
\begin{align}
|I_{\Gamma^{(2)}_{q,\xi}}(p/q,\ga)|\le & C\frac{e^{-\ga\xi}}{\xi}\le C\frac{e^{-\frac{p}{2q}\xi}}{\xi}
\rightarrow 0,\qquad\mbox{as $\xi\to\infty$},\\
|I_{\Gamma^{(4)}_{q,\xi}}(p/q,\ga)|\le & C\frac{e^{-(1+p/q-\ga)\xi}}{\xi}\le C\frac{e^{-\frac{p}{2q}\xi}}{\xi}
\rightarrow 0,\qquad\mbox{as $\xi\to\infty$}.
\end{align}
and we conclude that
\be\la{limI}
I(p/q,\ga)=\lim_{\xi\to\infty}I_{\Gamma_{q,\xi}}(p/q,\ga).
\ee
On the other hand, $I_{\Gamma_{q,\xi}}(p/q,\ga)$ is given by the sum of residues inside the
contour. Clearly, the integrand has poles there at the points:
\begin{align}
z_m &=2\pi i m,\qquad m=1,\dots,q-1\\
\wt z_n &=2\pi i {q\over p}(n-1/2),\qquad n=1,\dots,p.
\end{align}
Note that for all these $m,n$, $z_m\neq \wt z_n$ because $m\neq {q\over p}{2n-1\over 2}$ as $q$ is odd.
Hence we conclude that all the poles inside $\Gamma_{q,\xi}$ are simple.
Computing the residues and using the facts that 
\[
{q\over p}=a+{p'\over q'},\qquad q'=p,
\] 
we obtain (\ref{1intrec}) by (\ref{limI}).  Note that the conditions on $p'$, $q'$ in $T(p'/q')$ are fulfilled 
since $q$ is odd.

Now, in order to obtain the inequality (\ref{1intest}), we evaluate the integral along the contour
$\Gamma_q$. On the lower part of it,
\begin{align}
|I_{{\pi i\over 2}+\mathbb{R}}(p/q,\ga)|&\le
2\int_0^{\infty}\frac{e^{-\ga x}+e^{-(1+{p\over q}-\ga)x}}
{(1+2e^{-{p\over q}x}\cos\frac{\pi p}{2q}+e^{-2{p\over q}x})^{1/2}
(1+e^{-2x})^{1/2}}\frac{dx}{(x^2+\frac{\pi^2}{4})^{1/2}}\nonumber\\
&< 4\int_0^{\infty}\frac{e^{-{p\over 2q}x}}{(x^2+\frac{\pi^2}{4})^{1/2}}dx=
4\int_0^{\infty}\frac{e^{-u}}{(u^2+(\frac{\pi p}{4 q})^2)^{1/2}}dx.\la{esti1}
\end{align}
Separating the final integral into 2 parts, one along $(0,1)$ and another along $(1,\infty)$, 
we can continue (\ref{esti1}) as follows:
\begin{align}
&< 4\left(\int_0^1\frac{du}{(u^2+(\frac{\pi p}{4 q})^2)^{1/2}}+
\frac{1}{(1+(\frac{\pi p}{4 q})^2)^{1/2}}\int_1^\infty e^{-u}du\right)\nonumber\\
&= 4\left(-\ln\frac{\pi p}{4q}+\ln\left[1+\sqrt{1+\left(\frac{\pi p}{4 q}\right)^2}\right]+
\frac{e^{-1}}{(1+(\frac{\pi p}{4 q})^2)^{1/2}}\right)\nonumber\\
&<  4\left(\ln{q\over p}+\ln\left[{4\over\pi}+\sqrt{\left({4\over\pi}\right)^2+1}\right]+e^{-1}\right).\la{1intest1}
\end{align}

Similarly, we obtain (recall that $q>1$)
\begin{align}
|I_{2\pi i(q-{1\over 4})+\mathbb{R}}(p/q,\ga)|&\le
4\int_0^{\infty}\frac{e^{-u}}{(u^2+(\pi(p-\frac{p}{4 q}))^2)^{1/2}}dx\nonumber\\
&<\frac{4e^{-1}}{\pi(p-\frac{p}{4 q})}<\frac{4e^{-1}}{\pi p}\frac{1}{1-\frac{1}{4q}}<
\frac{5}{e\pi p}.\la{1intest2}
\end{align}
The sum of (\ref{1intest1}) and (\ref{1intest2})
gives
(\ref{1intest}), and thus we finish the proof of Lemma 5. $\Box$

We will need another similar lemma.

\medskip
\noindent
{\bf Lemma 6} {\it Let $(p,q)=1$, $p>0$, $q>1$,
\be\la{2int}
J(p/q,\ga,\de)=2\int_{\Gamma_q}\frac{(-1)^{\de}e^{(1+p/q)z}}
{(1-(-1)^{\de}e^{z p/q})(1+e^z)}\frac{e^{-\ga z}}{z}dz,
\qquad \frac{1}{2}\frac{p}{q}\le\ga\le1+\frac{1}{2}\frac{p}{q},
\ee
where $\de=\{0,1\}$, and $\Gamma_q$ is the same contour as in (\ref{1int}). 
Assume that $p$ is odd if $\de=0$, and either ($p$ -- even, $q$ -- odd), or ($p$ -- odd, $q$ -- even)
if $\de=1$.
Let $p/q$ have the continued fraction (\ref{cf}).
Then 
\be\la{2intrec}
J(p/q,\ga,\de)=T(p/q,\ga,\de)+
\begin{cases}
-S(p'/q',\ga'),& \mbox{if}\quad \de=0\\
(-1)^{\ep'}T(p'/q',\ga',a+1\,\mod 2),& \mbox{if}\quad \de=1
\end{cases},\qquad \ga'=\frac{q}{p}\ga (\mod 1),
\ee
where $\ep'=0$ if $\ga'=\frac{q}{p}\ga (\mod 2)$, and $\ep'=1$ otherwise.

Moreover, there holds the bound with $\beta$ from (\ref{1intest})
\be\la{2intest}
|J(p/q,\ga,\de)|<A_{\de}\left(4\ln\frac{q}{p}+\frac{5}{e\pi p}+
\beta\right),\qquad A_{\de}=
\begin{cases}
(1-\cos^2 \pi{p\over q})^{-1/2},& \mbox{if}\quad \de=0\\
1,& \mbox{if}\quad \de=1
\end{cases} 
\ee
}

\noindent
{\bf Proof}
We argue as in the proof of Lemma 5. The poles of the integrand in (\ref{2int})
inside ${\Gamma_{q,\xi}}(p/q,\ga)$ are:
\begin{align}
z_m &=2\pi i (m-1/2),\qquad m=1,\dots,q,\\
\wt z_n &=2\pi i {q\over p}(n-\de/2),
\end{align}
and $n=1,\dots,p-1$ if $\de=0$, while $n=1,\dots,p$ if $\de=1$. Our assumptions on the parity of
$p$, $q$ immediately imply that all $z_n\neq \wt z_m$ and hence all the poles inside  ${\Gamma_{q,\xi}}(p/q,\ga)$ are simple. Computing the residues we obtain (\ref{2intrec}).

Denote by $J_{\Gamma}(p/q,\ga,\de)$ the integral which has the same 
integrand as in (\ref{2int}), but with integration over some contour $\Gamma$.
As in the previous proof,
consider now an estimate for the integral along the lower part of $\Gamma_q$:
\begin{align}
|J_{{\pi i\over 2}+\mathbb{R}}(p/q,\ga)|&\le
2\int_0^{\infty}\frac{e^{-\ga x}+e^{-(1+{p\over q}-\ga)x}}
{(1-2(-1)^{\de}e^{-{p\over q}x}\cos\frac{\pi p}{2q}+e^{-2{p\over q}x})^{1/2}
(1+e^{-2x})^{1/2}}\frac{dx}{(x^2+\frac{\pi^2}{4})^{1/2}}\nonumber\\
&<4A_{\de}\int_0^{\infty}\frac{e^{-{p\over 2q}x}}{(x^2+\frac{\pi^2}{4})^{1/2}}dx,
\end{align}
where $A_{\de}$ is given in (\ref{2intest}). The rest of the argument is very similar to that in the proof of Lemma 5. $\Box$

\medskip
\noindent
{\bf Lemma 7} (recurrence) {\it 
Let $q>1$, $\frac{p}{q}=\frac{p_n}{q_n}=[a_1,\dots,a_n]$ and denote $t_j=[a_j,\dots,a_n]$, $j=1,\dots,n$. 
Assume that $a_1$ is odd, and $a_j$ are even for $j=2,\dots,n$. Fix $0\le k\le (q-1)/2$.
Then
\be\la{recsum}
S(p/q,\ga_1)=I(t_1,\ga_1)+\sum_{j=2}^n (-1)^{\ep_j} J(t_j,\ga_j,1)+2(-1)^{\ep_n}e^{-i\pi\ga_{n}a_n}
\ee
where $\ep_j=\{0,1\}$ and
\be
\ga_1=-2kt_1+k_1,\qquad \ga_j=k_{j-1}t_j+k_j,\quad j=2,\dots,n,
\ee
with the sequence $k_j\in\bbz$, $j=1,\dots, n$, chosen so that ${1\over 2}t_j\le\ga_j\le 1+{1\over 2}t_j$.

Furthermore, for this $p/q$, there holds the following bound for the sum (\ref{Sk2}):
\be\la{Sk2est}
|S_k|=\left|\sum_{m=1}^{q-1}\frac{1}{m}
\frac{\cos \pi{p\over q}m(4k+1)}{\cos \pi{p\over q}m}\right|
<\left( 4+\frac{\beta}{\ln 2}\right)\ln q+9.
\ee

If, in addition, $q_{k+1}\ge q_k^{\nu}$, for some $\nu>1$ and all $1\le k\le n-1$, then for any $\ve>0$
there exist $Q=Q(\ve,\nu)$ such that if $q>Q$, 
\be\la{Sk2estnew}
|S_k|<(4+\ve)\ln q.
\ee
}

\noindent
{\bf Proof}
First, note that $t_1=p/q$ and $t_j=1/(a_j+t_{j+1})$, $j=1,\dots,n-1$, $t_n=1/a_n$.
Now note a simple fact that the conditions $a_1$ -- odd, $a_2,\dots, a_n$ -- even ensure that
$q$ is odd and all the fractions $t_j$ for $j=2,\dots, n$ are either $\frac{\mbox{odd}}{\mbox{even}}$ or
$\frac{\mbox{even}}{\mbox{odd}}$.
We now choose $k_1$ so that $\ga_1=-2kt_1+k_1$ satisfies ${1\over 2}t_1\le\ga_1\le 1+{1\over 2}t_1$.
Applying Lemma 5 to $S(t_1,\ga_1)$, we obtain
\be\la{r1}
S(t_1,\ga_1)=I(t_1,\ga_1)+(-1)^{\ep}T(t_2,\ga_2, a_1\,\mod 2)
\ee
We can and will choose $\ga_2$ so that ${1\over 2}t_2\le\ga_2\le 1+{1\over 2}t_2$
by picking the appropriate $k'_2$ and setting 
\be
\ga_2={q\over p}\ga_1 +k'_2={q\over p}(-2kt_1+k_1)+k'_2=-2k+{q\over p}k_1+k'_2=
-2k+(a_1+t_2)k_1+k'_2=t_2k_1+k_2,  
\ee
where $k_2=-2k+a_1k_1+k'_2$. The constant $\ep$ is determined by $\ga_2$ as
described in Lemma 5.

Since $a_1\,\mod 2=1$, we can now apply Lemma 6 with $\de=1$ to $T$ in the r.h.s. of (\ref{r1}).
This gives
\be\la{r2}
S(t_1,\ga_1)=I(t_1,\ga_1)+(-1)^{\ep}[J(t_2,\ga_2,1)-(-1)^{\ep'}
T(t_3,\ga_3, a_2+1\,\mod 2)].
\ee
with 
\[
\ga_3=\frac{1}{t_2}\ga_2+k'_3=k_1+(a_3+t_3)k_2+k'_3=t_3k_2+k_3.
\]
chosen so that ${1\over 2}t_3\le\ga_3\le 1+{1\over 2}t_3$.

Note that according to our assumption $a_j+1\, \mbox{mod}\, 2 =1$, $j=2,\dots,n$, so that we can continue
applying Lemma 6 with $\de=1$ in (\ref{r2}) recursively. At the final step, revisiting the residue calculations
for Lemma 6 gives:
\be
T(t_n,\ga_n,1)=J(t_n,\ga_n,1)+2e^{-i\pi\ga_n a_n},
\ee
which proves (\ref{recsum}).

Now using the bounds (\ref{1intest}) and (\ref{2intest}) for $\de=1$, we write
\begin{align}\la{Sk2est-1}
|S_k|&<|S(p/q,\ga_1)|\le
|I(t_1,\ga_1)|+\sum_{j=2}^n |J(t_j,\ga_j,1)|+2<
4\ln\frac{1}{t_1t_2\cdots t_n}+\frac{5}{e\pi}\sum_{j=1}^n\frac{1}{p'_j}+n\beta+2,
\end{align}
where $t_j=[a_j,\dots,a_n]=\frac{p'_j}{q'_j}$.

Recall that we denote $q=q_n$, $p=p_n$.
Observe that the following recurrence \ci{HL1922} with $t_{n+1}=0$
\[
q_n+t_{n+1}q_{n-1}=(a_n+t_{n+1})q_{n-1}+q_{n-2}=\frac{1}{t_n}(q_{n-1}+t_n q_{n-2})
\]
gives
\be\la{t-identity}
\frac{1}{t_1t_2\cdots t_n}=q_n.
\ee
Note that this equation holds for any continued fraction.
 
Furthermore, the recurrence $t_{j-1}=1/(a_{j-1}+t_j)$ implies
\[
p'_{j-1}=q'_j,\qquad q'_{j-1}=a_{j-1}p'_{j-1}+p'_j,
\]
and so
\[
p'_{j-2}=a_{j-1}p'_{j-1}+p'_j,
\]
which, with the initial conditions $p'_n=1$, $p'_{n-1}=a_n$, 
and our assumption that all $a_j$, $j=2,\dots n$, are even, gives
\be
p'_j\ge a_{j+1}\cdots a_n\ge 2^{n-j},\qquad j=1,\dots,n.
\ee
Thus we have
\be\la{sp}
\sum_{j=1}^n\frac{1}{p'_j}<2.
\ee

Finally, since
\[
q_n=a_n q_{n-1}+q_{n-2}\ge a_n q_{n-1}\ge \cdots\ge a_na_{n-1}\cdots a_{1}\ge 2^{n-1},
\]
we have
\be\la{nineq}
n\le\frac{\ln q_n}{\ln 2}+1.
\ee
(In fact, as is well known, a slightly worse bound on $n$ holds for any continued fraction.)

Using (\ref{t-identity}), (\ref{sp}), (\ref{nineq}) in (\ref{Sk2est-1}), we 
obtain (\ref{Sk2est}). 

To obtain (\ref{Sk2estnew}) note first the following.

\medskip
\noindent{\bf Lemma 8} {\it
Let $\nu>1$, $n\ge 2$,
$p_k/q_k=[a_1,a_2,\dots,a_k]$, and such that $q_{k+1}\ge q_k^\nu$, $1\le k\le n-1$, $q_2\ge 3$.
Then
\[
n\le \frac{1}{\ln\nu}\ln\frac{\ln q_n}{\ln 3}+2
\]
}

\noindent
{\bf Proof} We have
\[
q_n\ge q_{n-1}^{\nu}\ge q_{n-2}^{\nu^2}\ge 
\cdots\ge q_2^{\nu^{n-2}}\ge 3^{\nu^{n-2}},\qquad n\ge 2,
\]
from which the result follows. $\Box$

Now using Lemma 8 instead of (\ref{nineq}) in (\ref{Sk2est-1}) and choosing $Q$ sufficiently large,
we obtain (\ref{Sk2estnew}) if $q>Q$, and finish the proof of Lemma 7. $\Box$

Bringing together (\ref{abs-si}), (\ref{Lend}), (\ref{1sumest}), and (\ref{Sk2est}),
yields the bound
\be
|\si'(0)|<e^{-{2\over 3}\ga_0}\frac{2}{3}q\cdot q^{\ga_0+4+\beta/\ln 2} e^{2\ga_0+9}
<q^{\ga_0+5+\beta/\ln 2}{2\over 3}e^{9+4\ga_0/3},\qquad q\ge 3.
\ee
Since $\ga_0+5+\beta/\ln 2=13.8\dots$ $<14$, 
${2\over 3}e^{9+4\ga_0/3}<e^{10}$,
and $\si(E)=-E$ if $q=1$, 
we obtain (\ref{si0}). 
Finally, using (\ref{Sk2estnew}) instead of (\ref{Sk2est}), we obtain (\ref{siw-new}).
This finishes the proof of Lemma 2. $\Box$

\section{Proof of Theorem 4}
Note first that since for any irrational $\al$
\[
\frac{1}{2q_nq_{n+1}}<\frac{1}{q_n(q_n+q_{n+1})}<\left|\al-{p_n\over q_n}\right|,
\]
we obtain that for any $\al$ satisfying the conditions of Theorem 4,
\be\la{qq}
\frac{C_3}{2}q_n^{\varkappa-1}<q_{n+1},\qquad n=1,2,\dots
\ee

We will denote $\mu_j(p/q)$, $w_j(p/q)$, etc, the values $\mu_j$, $w_j$, etc, for the spectrum $S(p/q)$.
Fix $n\ge 1$.
Let $E_2$, $E_0$ be the right edges of the centermost bands in $S(p_n/q_n)$, $S(p_{n+1}/q_{n+1})$,
respectively. By Lemma 2,
\be\la{E2}
E_2=\mu_0(p_{n}/q_{n})=w_0(p_{n}/q_{n})>\frac{4}{C_2 q_n^{C_1}}.
\ee
On the other hand by (\ref{w-est}) and (\ref{qq}), we have
\be\la{E0}
E_0=\mu_0(p_{n+1}/q_{n+1})=w_0(p_{n+1}/q_{n+1})<\frac{4e}{q_{n+1}}
<\frac{8e}{C_3 q_n^{\varkappa-1}}.
\ee
We will now show that $G_0(p_{n+1}/q_{n+1})\subset (E_0, E_2-\epsilon)$ for a suitably chosen $C_3$. 

Recall a continuity property found by Avron, Van Mouche, Simon \ci{AMS}: if $E\in S(\beta)$, there is
$E'\in S(\beta')$ such that
\be\label{continuity}
|E-E'|< C |\beta-\beta'|^{1/2}.
\ee
In \ci{AMS}, the authors give a good bound on $C$ requiring that $|\beta-\beta'|$ be sufficiently small.
As the reader can verify, a trivial modification of the proof in \ci{AMS} allows us to fix $C=60$ 
for the almost Mathieu operator
(worse than in \ci{AMS}) but without any condition on $\beta,\beta'\in (0,1)$. Thus we set $C=60$. 

This continuity property (\ref{continuity}) for $\beta=p_n/q_n$, $\beta'=p_{n+1}/q_{n+1}$,
together with the identity
\[
\left|\frac{p_n}{q_n}-\frac{p_{n+1}}{q_{n+1}}\right|=\frac{1}{q_nq_{n+1}}
\]
and the bound (\ref{qq})
implies that there exists $E'\in S(p_{n+1}/q_{n+1})$ such that
\be
E'\in \left({E_2\over 2}-{C\over \sqrt{q_nq_{n+1}}},{E_2\over 2}+{C\over \sqrt{q_nq_{n+1}}}\right)
\subset
\left({E_2\over 2}-{C\over \sqrt{C_3 q_n^{\varkappa}/2}},{E_2\over 2}+{C\over \sqrt{C_3 q_n^{\varkappa}/2}}\right).
\ee
Using (\ref{E2}) and recalling that\footnote{Note that here just $\varkappa=2C_1$ would do, cf a remark
following Theorem 4.}
$\varkappa=4C_1$, we see that 
\[
{E_2\over 2}-{C\over \sqrt{C_3 q_n^{\varkappa}/2}}>\frac{2}{C_2 q_n^{C_1}}-
{C\over \sqrt{C_3/2} q_n^{2C_1}}=\frac{2}{C_2 q_n^{C_1}}\left(1-{C C_2\over \sqrt{2C_3} q_n^{C_1}}\right),
\]
and setting now 
\be
C_3=4^2 C^2 C_2^4=4^2 60^2 C_2^4,
\ee
we have
\be\la{E*}
{E_2\over 2}-{C\over \sqrt{C_3 q_n^{\varkappa}/2}}>
\frac{q_n^{-C_1}}{C_2}. 
\ee
On the other hand, using (\ref{E0}), we obtain that
\[
E_0< \frac{8e}{C_3 q_n^{4C_1-1}}=\frac{e}{2C^2 C_2^4 q_n^{4C_1-1}}<
\frac{q_n^{-C_1}}{C_2}.
\]
Inequality (\ref{E*}) also shows that
\be\label{E2E}
{E_2\over 2}+{C\over \sqrt{C_3 q_n^{\varkappa}/2}}<E_2.
\ee
Thus,
\[
E'\in (E_0,E_2),
\]
which implies that
\[
G_0(p_{n+1}/q_{n+1})\subset (E_0, E_2-\epsilon),
\]
for some $\epsilon>0$. The corresponding result for $G_{-1}$ follows by the symmetry
of the spectra. This proves the statement (a) of Theorem 4.

Now by the continuity (\ref{continuity}) with $\beta=\al$, $\beta'=p_n/q_n$, and Theorem 3, 
we conclude that, for all $n=1,2,\dots$, there exists a gap $G_{n,2}(\al)$ of $S(\al)$ such that 
$G_{n,2}(\al)\cap G_0(p_n/q_n)\neq\emptyset$ and of length
\be
\De_{n,2}(\al)>\De_0(p_n/q_n)-2C|\al-p_n/q_n|^{1/2}>
\frac{1}{C_2^2 q_n^{2C_1}}-
\frac{2C}{C_3^{1/2}q_n^{\varkappa/2}}=
\frac{1}{2 C_2^2 q_n^{\varkappa/2}}.
\ee

We now verify that the gaps $G_{n,2}(\al)$, $G_{n+1,2}(\al)$ are distinct.
Using the continuity once again, we obtain that there exists a point $E''\in S(\al)$ such that
\be\label{Epp}
E''\in \left({E_2\over 2}-{C\over \sqrt{C_3 q_n^{\varkappa}}},{E_2\over 2}+
{C\over \sqrt{C_3 q_n^{\varkappa}}}\right).
\ee
Now it is easy to verify, similar to the calculations above, that 
\[
{E_2\over 2}-{C\over \sqrt{C_3 q_n^{\varkappa}}}>{7\over 4}{1\over C_2q_n^{C_1}}
\]
and
\[
E_0+\frac{2C}{\sqrt{C_3q_n^\varkappa}}<{1\over C_2q_n^{C_1}},
\]
and therefore (\ref{Epp}) together with (\ref{E2E}) yields
\be
E''\in \left(E_0+\frac{2C}{\sqrt{C_3q_{n+1}^\varkappa}},E_2\right).
\ee
Thus $G_{n+1,2}(\al)$ lies to the left of $E''$, and $G_{n,2}(\al)$ to the right of $E''$,
so that $G_{n,2}(\al)$ and $G_{n+1,2}(\al)$ are distinct gaps, $n=1,2,\dots$. 
Similar results for $G_{n,1}(\al)$ follow by the symmetry.
This proves the statement (b) of Theorem 4.

The proof of the statement (c) is similar and based on (\ref{siw-new}). It is a simple exercise. $\Box$

\section*{Acknowledgement}
The work of the author was partially supported by the Leverhulme Trust research fellowship
RF-2015-243. The author is grateful to Jean Downes and Ruedi Seiler for their hospitality at TU Berlin where part
of this work was written and to the referees for very useful comments.

\end{document}